\documentclass[letterpaper, 10 pt, conference]{ieeeconf}  % Comment this line out if you need a4paper

\IEEEoverridecommandlockouts % This command is only needed if  you want to use the \thanks command

\usepackage{afterpage}
\usepackage{epsfig} % for postscript graphics files
\usepackage{times} % assumes new font selection scheme installed
\usepackage{amsmath} % assumes amsmath package installed
\usepackage{amssymb}  % assumes amsmath package installed
\usepackage{amsfonts}
\usepackage{mathptmx} % assumes new font selection scheme installed
\usepackage{tabularx}
\usepackage{epstopdf}
\usepackage{float}
\usepackage{epic,color}
\usepackage{mathrsfs}
\usepackage{dsfont,MnSymbol}
\usepackage{algorithm}
\usepackage{multirow} % Multirow for tables
\usepackage[T1]{fontenc}

\newcounter{rmnum}

\newcounter{anum}

\def\IEEEQEDclosed{\mbox{\rule[0pt]{1.3ex}{1.3ex}}}

\def\qed{\ifmmode\IEEEQEDclosed\else{\unskip\nobreak\hfil
		\penalty50\hskip1em\null\nobreak\hfil\IEEEQEDclosed
		\parfillskip=0pt\finalhyphendemerits=0\endgraf}\fi}
\def\qed{\hspace*{\fill}~\IEEEQED\par\endtrivlist\unskip}%\mbox{\rule[0pt]{1.3ex}{1.3ex}}}

\def\Re{\mathbb{R}}

\def\Sec#1{Sec.~\ref{#1}}

\def\notes#1{\marginpar{\tiny #1}\typeout{Notes!
Notes!
Notes!
}}
\renewcommand{\notes}[1]{\typeout{notes!}}
\def\FRAC#1#2#3{\genfrac{}{}{}{#1}{#2}{#3}}
\def\half{{\mathchoice{\FRAC{1}{1}{2}}%
{\FRAC{2}{1}{2}}%
{\FRAC{3}{1}{2}}%
{\FRAC{4}{1}{2}}}}

\def\Re{\field{R}}

\def\Sec#1{Sec.~\ref{#1}}

\def\clP{{\cal P}}
\def\clZ{{\cal Z}}

\def\Sec#1{Sec~\ref{#1}}

\def\E{{\sf E}}

\def\Sec#1{Sec.~\ref{#1}}

\def\IEEEQEDclosed{\mbox{\rule[0pt]{1.3ex}{1.3ex}}}
\def\qed{\nobreak\hfill\IEEEQEDclosed}

\def\clZ{{\cal Z}}

\newtheorem{theorem}{Theorem}
\newtheorem{example}{Example}
\newtheorem{definition}{Definition}

\newtheorem{remark}{Remark}
\newtheorem{proposition}{Proposition}

\def\beq{\begin{eqnarray}} 
\def\bc{\begin{center}} 
\def\be{\begin{enumerate}}
\def\bi{\begin{itemize}} 
\def\bs{\begin{small}}
\def\bS{\begin{slide}}
\def\ec{\end{center}} 
\def\ee{\end{enumerate}}
\def\ei{\end{itemize}}
\def\es{\end{small}}
\def\eS{\end{slide}}
\def\eeq{\end{eqnarray}}

\newcommand{\newP}[1]{\medskip\noindent{\bf #1:}}

\newcommand{\ud}{\,\mathrm{d}}

\def\Re{\mathbb{R}}
\def\E{{\sf E}}

\def\Sec#1{Sec.~\ref{#1}}

\def\Prop#1{Prop.~\ref{#1}}

\def\clP{{\cal P}}
\def\clZ{{\cal Z}}

\renewcommand{\Re}{\mathbb{R}}

\def\FRAC#1#2#3{\genfrac{}{}{}{#1}{#2}{#3}}

\def\sJ{{\sf J}}
\def\sP{{\sf P}}

\def\ones{{\sf 1}}

\def\clF{{\cal F}}
\def\clG{{\cal G}}

\def\clP{{\cal P}}
\def\clU{{\cal U}}

\def\clZ{{\cal Z}}
\def\E{{\sf E}}

\def\bS{\mathbb{S}}

\newlength{\noteWidth}
\long\def\notes#1{\ifinner
	{\tiny #1}
	\else
	\marginpar{\parbox[t]{\noteWidth}{\raggedright\tiny #1}}
	\fi}

\title{\LARGE \bf The Conditional Poincar\'e Inequality for Filter Stability}

\author{Jin Won Kim, Prashant G. Mehta and Sean Meyn% <-this % stops a space
	\thanks{Financial support from the 
		NSF grant 1761622 and the 
		ARO grant W911NF1810334 is gratefully acknowledged. 
	}% <-this % stops a space
	\thanks{J.~W. Kim and P.~G.~Mehta are with the Coordinated
		Science Laboratory and the Department of Mechanical Science and
		Engineering at the University of Illinois at Urbana-Champaign
		(UIUC);  S.~Meyn is with the Department of
                Electrical and Computer Engineering at the University
                of Florida at Gainesville; Corresponding email: mehtapg@illinois.edu.}
}

\begin{document}

%Deterministic energy and variance, Conditional energy and variance
\def\den{\operatorname{enr}^{\bar{\mu}}}
\def\dvar{\operatorname{var}^{\bar{\mu}}}
\def\cen{{\cal E}^{\bar{\mu}}}
\def\cvar{{\cal V}^{\bar{\mu}}}

\maketitle
\thispagestyle{empty}
\pagestyle{empty}

%%%%%%%%%%%%%%%%%%%%%%%%%%%%%%%%%%%%%%%%%%%%%%%%%%%%%%%%%%%%%%%%%%%%%%%%%%%%%%%%
\begin{abstract}

This paper is concerned with the problem of nonlinear filter stability
of ergodic Markov processes.  The main contribution is the conditional
Poincar\'e inequality (PI), which is shown to yield filter stability.
The proof is based upon a recently discovered duality which is used to
transform the nonlinear filtering problem into a 
stochastic optimal control problem for a backward stochastic
differential equation (BSDE).  Based on these dual formalisms, a comparison is drawn between the
stochastic stability of a Markov process and the filter stability.
The latter
relies on the conditional PI described in this paper, whereas the
former relies on the standard form of PI.

% The former is based upon assuming the standard form of PI whereas the latter relies on the conditional PI introduced in this paper. 

\end{abstract}
%%%%%%%%%%%%%%%%%%%%%%%%%%%%%%%%%%%%%%%%%%%%%%%%%%%%

\section{Introduction}
\label{sec:intro}

%Literature survey.  

The Poincar\'e (or spectral gap) inequality (PI) is central to the
subject of stochastic stability of Markov
processes~\cite[Ch.~4]{bakry2013analysis}.  The PI is the simplest
condition which quantifies ergodicity and convergence to
stationarity: The Poincar\'e constant gives the rate of exponential
decay.  Apart from stochastic stability, the PI has a rich history.
It is the fundamental inequality in the study of the elliptic PDEs.

% During the past decade, spurred by advances in control-oriented
% algorithms for particle filtering (viz., the feedback
% particle filter (FPF)~\cite{taoyang_TAC12}), PI has become
% increasingly relevant in the numerical approximation of nonlinear
% filters.  In a numerical implementation, the FPF requires a solution
% of an elliptic PDE termed the Poisson equation of nonlinear filtering~\cite{variational_FPF}.  The PI is needed to ensure
% well-posedness of its solution~\cite{taghvaei2020diffusion}. Although the PI is simply assumed to
% hold in most applications of FPF, explicit bounds have been
% obtained in some special cases~\cite{variational_FPF,pathiraja2020mckean}.

%Rigorous results for the general case remain largely open.  

The goal of this paper is to propose a generalization of the PI for
the purpose of nonlinear filter stability analysis.  Specifically, a continuous-time filtering model is
considered where (a) the Markov process is ergodic, and (b)
the observations are corrupted by additive white noise.  
% ergodic Markov processes where the observations are corrupted by
% additive white noise (this is the setting also for the FPF).  
In the study of the Wonham filter, this model is referred to as the
{\em
  ergodic signal} case.  A companion paper, also published in the
proceedings of this conference, tackles the more general
non-ergodic signal case~\cite{kim2021detectable}.

For the case of ergodic signal, a pioneering early contribution
is~\cite{ocone1996asymptotic}.  Early work is based
on contraction analysis of the random matrix products arising from
recursive application of the Bayes' formula~\cite{atar1997lyapunov}
(see also~\cite[Ch. 4.3]{Moulines2006inference}).  Using related
techniques, the analysis of the Zakai equation leads to useful 
formulae for the Lyapunov exponents under assumptions on model
parameters and noise limits~\cite{atar1997exponential, Atar1999}.
%, and
%convergence rates estimated using Feynman-Kac type representation~\cite{Atar1999}.
% , and stability results for a special class of (Bene\v{s})
% filters was also obtained~\cite{ocone1999asymptotic}.  
For the ergodic signal case, a comprehensive account appears
in~\cite{budhiraja2003asymptotic} and the first complete solution
is given in~\cite{baxendale2004asymptotic}.  A necessary and
sufficient characterization of the asymptotic properties
of the Wonham filter, for both ergodic and non-ergodic cases,
appears in~\cite{van2009observability}.  For an accessible account of the
problem and the solution methods, see the review
papers~\cite{chigansky2006stability,chigansky2009intrinsic} and
references therein.

% A number of results have been obtained
% including bounds on the Lyapunov exponents of the Zakai equation (the
% un-normalized filter) in~\cite{atar1997exponential}, the formula for
% the relative entropy and its interpretation as a Lyapunov function for
% the filter in~\cite{clark1999relative},and  stability results on
% special class of (Bene\v{s}) filters in~\cite{ocone1999asymptotic}.
% For the ergodic signal case in particular, a comprehensive account of
% the asymptotic properties of the filter appears in~\cite{budhiraja2003asymptotic}.  The first complete
% solution of the ergodic 

% and application of
% the intrinsic methods for proving filter stability
% in~\cite{baxendale2004asymptotic}.  
% A complete characterization appears
% in~\cite{van2009observability} as a special case of more general
% results.  Some recent extensions of these prior 
% results appear in~\cite{mcdonald2018stability,mcdonald2019cdc}.   

Notwithstanding the fundamental importance of the PI for 
% {\em both} (a)
stochastic stability of Markov processes, 
% , and (b)
% numerical approximation of the nonlinear
% filter,
it is {\em not} a major theme in the filter stability
literature.  The closest is the appearance of the Brascamp-Lieb (B-L) inequality in~\cite[Chapter 4]{van2006filtering}
and~\cite{stannat2005stability}.  In both these references, the B-L
inequality is employed for the filter stability analysis of the It\^o
diffusions assuming linear observations.  The main point 
is that the B-L inequality is {\em not} central to the convergence
analysis.  In particular, the convergence rate relies on the uniform
convexity of the optimal value function in~\cite{van2006filtering} and
(equivalently) the log concavity of the posterior 
in~\cite{stannat2005stability} (and this is primarily a result of the
linear observation model). 
% (and this is primarily a result of the linear observation model).
The value function and the posterior are related through the log
transformation which is how duality in the nonlinear filters is
historically
understood~\cite{FlemingMitter82,mitter2003,todorov2008general}.
Apart from these works, it is known that in the large noise
limit, the top Lyapunov exponent converges to the spectral gap of the
ergodic Markov process~\cite{delyon1988lyapunov}.  
% provide a connection to the PI.  

% More importantly, unlike the PI, the use of the B-L inequality in these papers is not
% related to the rate of convergence.  

This paper has a single contribution: generalization of the PI for
Markov processes to the proposed conditional PI for the nonlinear filter.  The conditional PI is shown to play the same role for 
filter stability as the standard PI for stochastic stability.  The
proof technique relies on a recently discovered
duality result whereby the nonlinear filtering problem is cast as a
BSDE-constrained stochastic optimal control problem~\cite{kim2019duality}.  Using these
methods, we are able to derive {\em all} the prior results where
explicit convergence rates are obtained.  These are pointed in the paper.

The outline of the remainder of this paper is as follows: The problem
formulation appears in Sec.~\ref{sec:prelim}.  The PI and the conditional
PI are introduced in Sec.~\ref{sec:PI}.  The filter stability results
appear in \Sec{sec:filt_stab}. The Appendix contains the proofs.% Additional discussions as well as all the
% proofs appear in the Appendix.  

\section{Problem formulation}\label{sec:prelim}

\newP{Notation} The state-space $\bS:=\{1,2,\hdots\}$ is either finite
(with $d$ elements) or is countable.
The set of probability vectors on $\bS$ is denoted by $\clP(\bS)$:
$\mu\in \clP(\bS)$ if $\mu(x)\geq 0$ and $\sum_{x\in\bS} \mu(x) =
1$. The space of functions on $\bS$ is denoted $C(\bS)$ and the
space bounded functions on $\bS$ is denoted as $C_b(\bS)$: $f\in C_b(\bS)$ if $|f(x)|<C$ for some constant $C$ for all $x\in\bS$.   When the cardinality of
$\bS$ is finite ($d$), the space of functions on $\bS$ is identified with $\Re^d$.  For a measure $\mu \in {\cal P}(\mathbb{S})$ and a function $f\in C_b(\bS)$,
$\mu(f):=\sum_{x\in\bS} \mu(x) f(x)$. For two functions $f,h\in
C_b(\bS)$, $f \, h$ is the function obtained through element-wise
product: $(f\,h)(x) = f(x)h(x)$ for all $x\in\bS$ and similarly
$f^2=f f$ is the square of the function. The functions of all
ones is denoted as $\ones$, i.e., $\ones(x)=1$ for all $x\in\bS$.

\def\lsq#1#2{L_{#1}^2([0,T]\,;{#2})}

\subsection{Filtering model}

Consider a pair of continuous-time stochastic processes $(X,Z)$
defined on a probability space $(\Omega,\clF,{\sf P})$. 
The state process $X:=\{X_t\in \mathbb{S}:t\ge 0\}$ is a stationary Markov
process with generator $A$ and an everywhere positive invariant measure $\bar\mu
\in \clP(\bS)$: $\bar\mu(x)>0$ for all
$x\in\bS$ and 
$\bar\mu(Af)=0$ for all $f\in C_b(\bS)$ (this is the case for any irreducible positive recurrent
Markov process).  The observation process $Z = \{Z_t\in\Re^m:t\geq 0\}$
is defined according to the following model:
\begin{equation}\label{eq:obs}
Z_t := \int_0^t h(X_s) \ud s + W_t
\end{equation}
where $h:\bS\to \Re^m$ is the observation function and 
$W=\{W_t\in\Re^m:t\geq 0\}$ is a Wiener process (w.p.) that is assumed
to be independent of $X$.  The covariance of $W$ is denoted $R$ which is assumed to be strictly positive-definite.  
% \footnote{The scalar-valued observation model
  % is considered for notation ease.}.
The filtration generated by $Z$
is denoted $\clZ := \{\clZ_t : 0\le t\le T\}$ where $\clZ_t =
\sigma(\{Z_s: 0\le s\le t\})$.

\newP{Function spaces} 
%Since $X$ is stationary, the
To stress the choice of the initial stationary prior $\bar\mu$, we
write the probability measure $\sP$ as $\sP^{\bar\mu}$, the
expectation is denoted ${\sf E}^{\bar\mu}(\cdot)$, and
$L^p(\bar\mu)$ is the space of random variables $Y$ with ${\sf
  E}^{\bar\mu}(|Y|^p)<\infty$.  The
space of square-integrable deterministic functions on
$\bS$ is denoted as $L^2(\bS)$: $f\in
L^2(\bS)$ if ${\sf
  E}^{\bar\mu}(|f(X_T)|^2) = \sum_x |f(x)|^2 \bar\mu(x) < \infty$.  The
space of square-integrable $\clZ_T$-measurable random
functions on $\bS$ is denoted as $L^2_{\clZ_T}(\bS)$:  $F\in
L^2_{\clZ_T} (\bS)$ if $F$ is $\clZ_T$-measurable and ${\sf
  E}^{\bar\mu}(|F(X_T)|^2) < \infty$.  Likewise, the space of
$\clZ$-adapted square-integrable S-valued stochastic processes is
denoted $L^2_{\clZ}([0,T];S)$.  Examples are $S=\Re^m$ for
vector-valued and $S=C(\bS)$ for function-valued stochastic
processes. 
% This notation is used not only for
% real-valued but also for vector
% and function-valued stochastic processes.

% Likewise, $\lsq{\clZ}{S}$
% denotes the space of $\clZ$-adapted square-integrable processes taking
% values in a measurable space $S$.  Common examples of $S$ are as
% follows: (i) $S=L^2(\bar\mu)$ in
% which case $\lsq{\clZ}{S}$ is the space of $\clZ$-adapted
% square-integrable real-valued stochastic
% processes (this space is more simply denoted as $L^2_{\clZ}([0,T])$); and (ii)  $S=L^2_{\clZ_T}(\bS)$  in which case
% $\lsq{\clZ}{S}$ is the space of $\clZ$-adapted square-integrable
% function-valued stochastic processes.  

% is the square-integrable
% random variables taking values in $S$. 

% For the filtration $\clZ$
% and a measurable space $S$, 
% Likewise, $L^2_{\clZ_T}(S)$ is the $\clZ_T$-measurable square-integrable
% random variables taking values in $S$.  

\medskip

The filtering problem is to compute the conditional 
distribution (posterior) of the state $X_t$ given $\clZ_t$.  The
posterior distribution at time $t$ is denoted $\pi_t^{\bar\mu}\in\clP(\bS)$.  For $f\in  C_b(\bS)$
\[
\pi_t^{\bar\mu}(f) := \E^{\bar\mu}(f(X_t)|\clZ_t)
\]
% The superscript is used to highlight the fact that the chain is
% stationary, i.e., the prior is the invariant measure
% $\bar\mu$.  

\subsection{Definition of filter stability}

For each $f\in  C_b(S)$, the Wonham filter is given by the stochastic
differential equation:
%~\cite{wonham1964}
\begin{equation}\label{eq:Wonham}
\ud \pi_t(f) = \pi_t(Af) \ud t + (\pi_t(hf) - \pi_t(h) \pi_t(f))R^{-1}(\ud Z_t- \pi_t(h) \ud t)
\end{equation}
With an initialization $\pi_0=\mu  \in \clP(\bS)$, the solution of the
Wonham filter is denoted as $\pi^\mu := \{\pi_t^\mu \in\clP(\bS):t\geq 0\}$.  The
posterior $\pi^{\bar\mu}$ results from the choice of
the initial condition $\pi_0 = \bar\mu$.  

\medskip

\begin{remark}\label{rem:Pmu_Pnu}
(see also the discussion in~\cite[Sec.~1]{chigansky2009intrinsic}.) Suppose
$\mu\in \clP(\bS)$ then (because $\bar\mu$ is everywhere positive)
$\mu \ll \bar\mu$ and setting the Radon-Nikodym (R-N) derivative
\[
\frac{\ud \sP^\mu}{\ud \sP^{\bar\mu}}(\omega) = \sum_{x} \frac{\mu(x)}{\bar\mu(x)}
\ones_{[X_0 = x]}(\omega)
\]
yields a new probability measure $\sP^\mu$ on the common measurable
space $(\Omega,\clF)$.  With respect to ${\sP}^\mu$, the expectation
operator is denoted ${\sf E}^{\mu}(\cdot)$ and
$L^p(\mu)$ is the space of random variables $Y$ with ${\sf
  E}^{\mu}(|Y|^p)<\infty$.  The solution of the Wonham filter
$\pi_t^{\mu}(f) = \E^{\mu}(f(X_t)|\clZ_t)$.  
\end{remark}

\medskip

\begin{definition}
The Wonham filter is {\em stable} if for each $f \in C_b(\bS)$, $\pi_T^\mu(f) \stackrel{L^1(\bar\mu)}{\longrightarrow}
\pi_T^{\bar\mu}(f)$, i.e., 
\begin{equation}\label{eq:filter-stability}
\E^{\bar\mu} \big(| \pi_T^\mu(f) - \pi_T^{\bar\mu}(f) |\big)
\longrightarrow 0 \quad \text{as}\; T\to\infty
\end{equation}
for all $\mu \in \clP_0 \subset \clP(\bS)$.  
\end{definition}

\medskip

Our goals are as follows: (i) characterize the subset $\clP_0 \subset
\clP(\bS)$ for which the filter stability
holds (ideally $\clP_0 = \clP(\bS)$); and (ii) obtain explicit estimates on the rate of convergence.  

\medskip

\begin{remark} Suppose $f\in C_b(\bS)$, $\mu,\nu \in \clP_0$ such
  that~\eqref{eq:filter-stability} holds.  Then
\begin{enumerate}
\item
Because $\pi_T^\mu(f) \stackrel{L^1(\bar\mu)}{\longrightarrow}
\pi_T^{\bar\mu}(f)$ and $\pi_T^\nu(f) \stackrel{L^1(\bar\mu)}{\longrightarrow}
\pi_T^{\bar\mu}(f)$, by the use of triangle
inequality, 
% \[
% \E^{\bar\mu}(|\pi_T^\mu(f) -  \pi_T^{\nu}(f)|) \leq
% \E^{\bar\mu}(|\pi_T^\mu(f) -  \pi_T^{\bar\mu}(f)|) + \E^{\bar\mu}(|\pi_T^\nu(f) -  \pi_T^{\bar\mu}(f)|) 
% \]
it also follows that $\pi_T^\mu(f) \stackrel{L^1(\bar\mu)}{\longrightarrow}
\pi_T^{\nu}(f)$. 
\item  (see~\cite[Remark 3.3]{ocone1996asymptotic}) Also $\pi_T^\mu(f) \stackrel{L^p(\bar\mu)}{\longrightarrow}
\pi_T^{\nu}(f)$ for all $p\geq 1$.  This is because, e.g., with
$p=2$,
\[
\E^{\bar\mu}(|\pi_T^\mu(f) -  \pi_T^{\nu}(f)|^2) \leq
2 \|f\|_{\infty} \, \E^{\bar\mu}(|\pi_T^\mu(f) -  \pi_T^{\nu}(f)|)  
\]
\item If $\pi_T^\mu(f) \stackrel{L^1(\bar\mu)}{\longrightarrow}
\pi_T^{\nu}(f)$ then (owing to Remark~\ref{rem:Pmu_Pnu}) it also
follows $\pi_T^\mu(f) \stackrel{L^1(\nu)}{\longrightarrow}
\pi_T^{\nu}(f)$ and therefore, $\pi_T^\mu(f) \stackrel{L^p(\nu)}{\longrightarrow}
\pi_T^{\nu}(f)$ for all $p\geq 1$.  
% \[
% \E^{\mu}(|\pi_T^\mu(f) -  \pi_T^{\nu}(f)|) \to 0
% \]
% This is because every $\mu\in\clP(\bS)$ is absolutely continuous with
% respect to $\bar\mu$.  
\end{enumerate}
\end{remark}

\medskip

\begin{remark}
The problem of stochastic stability is a special case when the
filtration $\clZ$ is independent of $X$  (e.g., $\clZ$ is trivial).  The counterpart of the Wonham filter is
the Kolmogorov's forward equation  
\begin{equation}\label{eq:Kol}
\ud \pi_t(f)  = \pi_t(Af) 
\end{equation}
With an initialization $\pi_0=\mu  \in \clP(\bS)$, the solution $\pi^\mu
:= \{\pi_t^\mu \in\clP(\bS):t\geq 0\}$ is now a deterministic process
(which serves to simplify the problem considerably).
Adapting~\eqref{eq:filter-stability} to this case, the basic problem of
stochastic stability is to show $\pi_T^\mu(f) \to \bar\mu(f)$ for all
$f\in C_b(\bS)$. 
\end{remark}

\section{Main assumption: Poincar\'e inequality (PI)}\label{sec:PI}

\subsection{Standard form of PI}

The carr\'e du champ operator is defined 
%for functions $f\in C_b(\bS)$
as
\[
\Gamma (f) (x) := \sum_{j \in \mathbb{S}} A(x,j) (f(x) - f(j))^2, \;\;
x \in \bS
\]
For a deterministic function $f\in L^2(\bS)$, the energy and variance are defined as
\begin{align*}
\den_0 (f) &:= \E^{\bar\mu} (\Gamma (f) (X_T)) = \sum_x  \bar{\mu}(x) \Gamma (f) (x)\\
\dvar_0(f) &:= \E^{\bar\mu} (|f(X_T) - \bar\mu(f)|^2) = \sum_{x} \bar{\mu}(x) |f(x) - \bar{\mu}(f)|^2 
\end{align*}
The standard form of the PI relates the two as follows:
\[
\text{PI}(\bar\mu,A):  \quad \den_0 (f) \geq c_0 \;
\dvar_0(f)\quad \forall\; f \in L^2(\bS)
\]
where $c_0>0$.  The PI is a standard assumption in the theory of Markov
processes to show stochastic stability~\cite{bakry2008rate}.

\subsection{Conditional form of PI}

For a random function $F\in L^2_{\clZ_T} (\bS)$, the energy and
variance are defined as follows:
\begin{flalign*}
& \text{(energy)} \quad \den_T(F)  := \E^{\bar\mu} (\Gamma (F) (X_T)) 
%= \E^{\bar\mu} ( \sum_{j \in \mathbb{S}} A(X_T,j) (F(X_T) - F(j))^2 )
\\
& \text{(variance)} \quad \dvar_T(F)  :=
\E^{\bar\mu} (|F(X_T) - \pi_T^{\bar\mu}(F)|^2) 
\end{flalign*}
%where $\pi_T^{\bar\mu}(F) = \E^{\bar\mu} (F(X_T)|\clZ_T)$. 

\begin{definition}\label{def:cPI}
A Markov process is said to satisfy the {\em conditional Poincar\'e
  inequality} (PI) 
%with respect to the filtration $\clZ$ and 
with a constant $c$ if 
\[
\text{PI}(\bar\mu,A; \clZ): \quad \den_T (F) \geq c \;
\dvar_T(F)\quad \forall \; F\in L^2_{\clZ_T}
(\bS),\;\forall\;T\geq 0
\]
\end{definition}

% \medskip

% \begin{remark} Recall the classical PI for a Markov chain:
% \[
% \sum_x \bar{\mu}(x) \Gamma (f) (x) \geq c \sum_{x} \bar{\mu}(x) |f(x)
% - \bar{\mu}(f)|^2 \;\; \forall \; f\in C_b(\bS)
% \]
% \begin{enumerate}
% \item Suppose $\clZ_T$ is trivial.  Then the conditional PI reduces to
%   the PI for a Markov chain.  
% \item Suppose the classical PI holds for a Markov chain with constant $c$.  Then conditional PI holds for all
%   deterministic choices of $F=f$:
% \begin{align*}
% {\cal E}^{\bar{\mu}} (f)  &=
%   \E^{\bar\mu} ( \sum_{j \in \mathbb{S}} A(X_T,j) (f(X_T) - f(j))^2 )\\
%                            & \geq  c \; \E^{\bar\mu} (|f(X_T) -
%                            \bar\mu(f)|^2)  \geq c \; \text{var}^{\bar{\mu}}(f) 
% \end{align*}
% \end{enumerate}
% \end{remark}
\medskip

\subsection{Examples}\label{ssec:examples}

%A simple example where the conditional PI holds.
\begin{example}\label{ex:2-states}
Suppose
$
A = \begin{bmatrix} -\lambda_1 & \lambda_1 \\ \lambda_2 & - \lambda_2 \end{bmatrix}
$
is irreducible.  Then $\lambda_1>0$ and $\lambda_2> 0 $. Observe that 
\begin{align*}
& \lambda_1 \pi_T^{\bar\mu}(1) (F(1) - F(2))^2 + 
\lambda_2 \pi_T^{\bar\mu}(2) (F(1) - F(2))^2 \\
& \geq
(\lambda_1+\lambda_2)  
\big(
\pi_T^{\bar\mu}(1) (F(1) - \pi_T^{\bar\mu}(F))^2  + \pi_T^{\bar\mu}(2) (F(2) - \pi_T^{\bar\mu}(F))^2 
\big)
\end{align*}
and therefore upon taking expectations on both sides,
\[
\den_T(F) \geq (\lambda_1+\lambda_2)  \dvar_T(F)
\]
Hence, the conditional PI holds for every irreducible 2-state Markov
chain with a constant $c=(\lambda_1+\lambda_2)$.  
Note that the best constant for standard PI is
$c_0=2(\lambda_1+\lambda_2)$. 

% It is noted that the constant $c=(\lambda_1+\lambda_2)$ is also tight,
% and in fact the same for both the standard and the conditional PI.
% Therefore, the conditional PI holds for every irreducible 2-state Markov
% chain.  
\end{example} 

\medskip

This simple example admits the following generalizations described in
the proposition with proofs in the Appendix~\ref{appdx:CPI-examples-pf}.  

\medskip

\begin{proposition}\label{prop:CPI-examples}
The conditional PI holds with the following constants (provided these are
positive)
\[
c = \sum_{j} \min_{i\in\bS:\;\; i\neq j} A(i,j), \quad c = \min_{i\neq j} \sqrt{A(i,j)}\sqrt{A(j,i)}
\]
\end{proposition}

\medskip

The first of the two formulae in \Prop{prop:CPI-examples} is related to the Doeblin type strong
mixing condition~\cite[Assumption 4.3.24]{Moulines2006inference}.   
The second formula is the same
as~\cite[Theorem 6]{atar1997lyapunov},~\cite[Theorem 4.3]{baxendale2004asymptotic} and~\cite[Corollary
2.3.2]{van2006filtering}.   
   
In this paper, the constant $c$ of the conditional PI is shown to play
the same role for the filter as the constant $c$ of the standard PI
does for the Markov process.  (Both give the rate of convergence.)  Moreover,
stronger results are possible with weaker notions of the conditional
PI (see Prop.~\ref{prop:beta_ind_ineq} and Example~\ref{ex:min_row}).

 % We have not seen the first
% of the two formulae described in the filter stability literature.  In
% stochastic stability, this constant is known for Doeblin chains
% (see Example~\ref{ex:exdoeblin}).    

% \medskip

% \begin{example}
% %\normalfont
% A Markov process is {\em Doeblin} if there exist a state $j^*\in \bS$
% such that $A(x,j^*) > c$ for all $x\in \bS \setminus \{j^*\}$.
% Then 
% % for
% % the Markov chain
% % \[
% % {\cal E}^{\bar{\mu}} (f) = \sum_i  \bar\mu(i) \sum_{j \in \mathbb{S}} A(i,j)
% % (f(i) - f(j))^2 \geq \epsilon \sum_i  \bar\mu(i) (f(i) - f(j^*))^2
% % \geq \epsilon \; \text{var}^{\bar{\mu}}(f)
% % \]
% % and similarly for the filter
% \begin{align*}
% {\cal E}_T^{\bar{\mu}} (F) & = 
%   \E^{\bar\mu} ( \sum_{j \in \mathbb{S}} A(X_T,j) (F(X_T) - F(j))^2 )\\
%   &\geq c\,
%   \E^{\bar\mu} ( (F(X_T) - F(j^*))^2 ) 
% %\\ 
% %& \geq \epsilon \; \E^{\bar\mu} ( (F(X_T) - \pi_T^{\bar\mu}(F) )^2 )  
% \geq c \, \text{var}_T^{\bar{\mu}}(F)
% \end{align*}
% Therefore, the conditional PI holds for Doeblin chains.  
% \end{example}

% \medskip

%\subsection{Relationship between standard and conditional PI}

\medskip

\begin{remark}
%Based upon results in the filter stability literature, 
One may
conjecture that the filter ``inherits'' the PI from the underlying
Markov process.  Note that the definition~\ref{def:cPI} is stated
for a general class of filtrations (not necessarily defined according
to the model~\eqref{eq:obs}).    
In the general settings, the conditional PI holds for
deterministic functions $f\in L^2(\bS)$.  This is because
%the conditional PI is inherited
%for the subspace $L^2(\bS)\subset L^2_{\clZ_T}(\bS)$ 
\begin{align*}
\den_T(f)  \geq  c_0 \; \E^{\bar\mu} (|f(X_T) -
                           \bar\mu(f)|^2)  \geq c_0 \; \dvar_T(f) 
\end{align*}
This shows that the PI and also the constant $c$ is inherited on
the subspace $L^2(\bS)\subset L^2_{\clZ_T}(\bS)$ of deterministic
functions.  
However, with general types of filtrations, it may not hold for random functions.  A counterexample
appears in the Appendix~\ref{apdx:counterexample}.
\end{remark}

% \medskip

% \begin{conjecture}
% Suppose the standard PI holds, with a constant $c$, for the Markov
% process $X$ and suppose the filtration $\clZ$ is obtained using the
% model~\eqref{eq:obs}.
% Then conditional PI holds also with constant $c$.  
% \end{conjecture}

\medskip

In the following, a proof of the filter stability is presented
based on the conditional PI.  
The proof of stochastic stability arises as a special case, and
is included in Appendix~\ref{sec:stoch_stab}.  The two sections are
self-contained and may be read independently of the other.  However, a
reader may benefit from reading the two sections simultaneously.  % The
% stochastic stability proof is a special case of the filter stability
% proof.     

% \subsection{Outline of the remainder of this paper}

% The outline of the remainder of this paper is as follows:
% \begin{enumerate}
% \item In \Sec{sec:stoch_stab}, a proof for the stochastic stability of
%   Markov processes is described based upon assuming the standard PI.  
% Although the result is also proved using direct means~\cite{}, our
%   proof relies on the use of a dual process. 
% \item In \Sec{sec:filt_stab}, a proof for filter stability is
%   described based upon assuming conditional PI.  The proof again is
%   based upon the use of duality.  
% \end{enumerate}
% A reader only interested in the filter stability result may directly
% skip ahead to \Sec{sec:filt_stab}.  Including \Sec{sec:stoch_stab} is
% useful for relating the PI and the conditional PI.  The stochastic
% stability proof is a special case of the filter stability proof.  

\medskip

\section{Filter Stability}
\label{sec:filt_stab}

For a random function $F\in L^2_{\clZ_T} (\bS)$, the conditional
energy and the conditional 
variance are defined as follows:
\begin{flalign*}
& \text{(cond. energy)} \quad \cen_T(F)  := \E^{\bar\mu} (\Gamma (F)
(X_T)|\clZ_T)  
%= \E^{\bar\mu} ( \sum_{j \in \mathbb{S}} A(X_T,j) (F(X_T) - F(j))^2 )
\\
& \text{(cond. var.)}\qquad \cvar_T(F)  := \E^{\bar\mu} (|F(X_T) - \pi_T^{\bar\mu}(F)|^2|\clZ_T)
\end{flalign*}
%By tower property of conditional expectation, these are related with
%the energy and the variance by
Upon taking expectations
\[
\den_T(F) = \E^{\bar{\mu}}(\cen_T(F)),\quad \dvar_T(F) = \E^{\bar{\mu}}(\cvar_T(F))
\]

\subsection{Duality}

In our prior work~\cite{peng1993backward}, the following backward stochastic differential
equation (BSDE) constrained optimal control problem is introduced.
The significance of the problem is that it is a dual of the nonlinear
filtering problem.

\newP{Dual optimal control problem}  
\begin{subequations}\label{eq:opt-cont-euclide}
	\begin{align}
          & \mathop{\text{Min}}_{U\in\;\clU}\;\; \sJ_T^{\bar\mu} (U)  =
            |Y_0(X_0)-\bar{\mu}(Y_0)|^2 + \E^{\bar\mu} \Big(\int_0^T \ell (Y_t,V_t,U_t\,;X_t) \ud t \Big)\label{eq:opt-cont-euclide-a}\\
          & \text{Subj.}\ - \ud Y_t(x)     = \big((A Y_t)(x) + h (x) (U_t +
            V_t(x))\big)\ud t - V_t^\top(x)\ud Z_t \nonumber \\
          & \quad\quad\quad\;\;\;  Y_T (x) = F(x) \;\; \forall \; x \in \bS \quad\text{(given)}
% \nonumber \\
	% &\quad\quad\;\;\;\; Y_T (x) \;\text{given}
\label{eq:opt-cont-euclide-b}
	\end{align}
\end{subequations}
where $\ell(y,v,u;x):= \Gamma(y)(x) + |u+v(x)|_R^2$, ${\cal U} :=
L^2_{\clZ}([0,T];\Re^m)$, and $F\in L^2_{\clZ_T}(\bS)$.  

\medskip

The existence and uniqueness of the optimal control follows from the
standard results in the BSDE constrained optimal control
theory~\cite{peng1993backward}.  The solution, including the formula
for optimal control, is described in~\cite[Theorem~1]{kim2019duality}. 
Let $U^{\text{opt}} :=\{U_t^{\text{opt}}\in\Re^m:\, 0\le t \le T\}$ be the
optimal control input and $({Y},{V}):=\{({Y}_t,{V}_t)\in C(\bS)\times C(\bS)^m:\,
0\le t \le T\}$ be the associated $\clZ$-adapted (optimal) trajectory
(the solution of the BSDE $(Y,V)\in L^2_{\clZ}([0,T]; C(\bS)\times C(\bS)^m)$).  Then 

%appears in the Appendix~\ref{apdx:dual-optimal-control}.

\medskip

%\begin{enumerate}
%\item
\noindent \textbf{1.}~\cite[Theorem~2]{kim2019duality}: The conditional mean 
\[
\pi_T^{\bar\mu}(Y_T) = \bar{\mu}({Y}_0) - \int_0^T U_t^{\text{opt}} \ud Z_t
\]
%Compare this with the Markov chain case where $\pi_T^{\bar\mu}(y_T)
%=\bar{\mu}(y_0)$.
  
\noindent \textbf{2.}~\cite[Theorem~5]{kim2019duality}:  Define a
$\clZ$-adapted process $M:=\{M_t:0\leq t\leq T\}$ as follows:
\begin{equation*}\label{eq:martingale}
M_t := \cvar_t(Y_t) - \int_0^t \big( \cen_s(Y_s) +
\sum_x \pi_s^{\bar\mu} (x) |U_s^{\text{opt}}+V_s(x)|^2_R \big) \ud s
%\big(\ell (Y_s,V_s,U_s^{\text{opt}}\,;X_s)\big) \ud s
\end{equation*}
Then $M$ is a $\sP^{\bar{\mu}}$-martingale.  
%With a different control
%$U\neq U^{\text{opt}}$, it is a super-mg.  
%\item~The variance of $Y_T$ is the optimal value function:

\medskip

\noindent \textbf{3.} Therefore, $\E^{\bar{\mu}}(M_T)
= \E^{\bar{\mu}}(M_0)$, which is expressed as
%expressed as 
% \[
% \text{var}_0^{\bar{\mu}}({Y}_0) + \E^{\bar\mu} \left( \int_0^T \ell (Y_t,V_t,U_t^{\text{opt}}\,;X_t)  \ud t\right)  = \text{var}_T^{\bar{\mu}}(Y_T)
% \] 
% %Again this equation has a direct parallel for the Markov chain.  
% %\end{enumerate}
% With the form of the cost function
\begin{align*}
% \text{var}_0^{\bar{\mu}}({Y}_0) + \E^{\bar\mu} \left( \int_0^T \Gamma
%   ({Y}_t) (X_t)  + |{U}_t^{\text{opt}} +{V_t}(X_t)|^2 \ud t\right)  =
%   \text{var}_T^{\bar{\mu}}({Y}_T) \\
% \therefore, &
\dvar_0({Y}_0) & + \int_0^T \den_t(Y_t) + \E^{\bar\mu}(|{U}_t^{\text{opt}}
  +{V_t}(X_t)|_R^2) \ud t  =
  \dvar_T({Y}_T)\\
& \therefore, \quad \dvar_0({Y}_0) + \int_0^T \den_t(Y_t) \ud t  \leq 
  \dvar_T({Y}_T)
\end{align*}
and using the conditional PI
\begin{equation}\label{eq:PI_ind_ineq_filter}
\dvar_0({Y}_0) \leq e^{-cT} \dvar_T({Y}_T)
\end{equation}
Equation~\eqref{eq:PI_ind_ineq_filter} is the backward inequality for
the variance of the dual process 
(this is the only place where the conditional PI is used).  The
counterpart for a Markov process is the
inequality~\eqref{eq:PI_ind_ineq_markov} for the dual process in Appendix~\ref{sec:stoch_stab}. %  As in
% \Sec{sec:stoch_stab}, this inequality is the only place in the stability
% proof where the PI is used.  

A more general result, described in the following proposition, is obtained
by considering the martingale $M$ directly.  Its proof appears in
Appendix~\ref{appdx:beta_ind_ineq_pf}.

\medskip

\begin{proposition}\label{prop:beta_ind_ineq}
Suppose $\beta=\{\beta_t: t \ge 0\}$ is any non-negative
$\clZ$-adapted process such that
\begin{equation}\label{eq:cPI-pathwise}
\cen_t(f) \ge \beta_t \; \cvar_t(f) \quad \sP^{\bar{\mu}}-\text{a.s.}
\quad \forall \;f \in L^2(\bS),\;\;
%\, F\in L^2_{\clZ_t}(\bS),\;
t\ge 0
\end{equation}
Then the backward inequality is of the form
\begin{equation}\label{eq:PI_ind_ineq_filter_stronger}
\dvar_0({Y}_0) \leq
\E^{\bar{\mu}}\left(e^{-\int_0^T \beta_t \ud t} \; \cvar_T(Y_T)\right) 
%\E^{\bar{\mu}}\left(\exp\Big(-\int_0^T \beta_t \ud t\Big) \; \cvar_T(Y_T)\right) 
\end{equation}
\end{proposition}

\medskip

Consequently, if $\frac{1}{T}\int_0^T \beta_t \ud t \to c$ (a
deterministic constant), then the backward
inequality~\eqref{eq:PI_ind_ineq_filter} for the variance is obtained
asymptotically.   
The following example shows how to choose $\beta$ to obtain an 
asymptotic formula for the convergence rate. 
%Returning to the question of choosing $\beta$, the
%following example shows that choosing $\beta$ as a stochastic process is not only possible but can also give a useful asymptotic formula for the convergence rate.  

\medskip

\begin{example}\label{ex:min_row}
	It is a straightforward calculation
	to verify
	\[
	\cen_t(f) \geq \left(\sum_{i} \pi_t^{\bar\mu}(i) \min_{j\in\bS:\;\; i\neq j}
	A(i,j) \right) \; \cvar_t(f),\quad \forall\;t\geq 0
	\]
	Set $\beta_t = \sum_{i} \pi_t^{\bar\mu}(i) \min_{j:i\neq
		j}A(i,j) $.  Using~\cite[Eq.~(5.15)]{baxendale2004asymptotic}, it is
	known that
	\[
	\lim_{T\to \infty} \frac{1}{T}\int_0^T \beta_t \ud t =  \sum_{i}
	\bar{\mu}(i) \min_{j\in\bS:\;i\neq j}A(i,j)  \quad \text{a.s.}
	\]
	The righthand-side then gives the asymptotic constant $c$ for the variance
	inequality~\eqref{eq:PI_ind_ineq_filter}.  This formula for the
	asymptotic convergence rate of the Wonham filter can be found 
	in~\cite[Theorem 4.2]{baxendale2004asymptotic}.  
	
\end{example}

\medskip
\begin{remark}
Inequality~\eqref{eq:cPI-pathwise} is the pathwise version of the
conditional PI. Note the inequality needs to be specified only for
deterministic functions.  Because of its pathwise nature, the
inequality then also holds for random functions $F\in  L^2_{\clZ_t}
(\bS)$.  A formal definition of pathwise PI is stated next.
%because
%\begin{align*}
%\cen_t(F) & = \sum_x  \pi_t^{\bar{\mu}}(x) \Gamma (F) (x) \\
%\cvar_t(F) &= \sum_{x} \pi_t^{\bar{\mu}}(x)  |F(x) - \pi_t^{\bar{\mu}}(F)|^2 
%% \begin{align*}
%% \cen_t(F) &= \sum_{i,j} \pi_t^{\bar\mu}(i) A(i,j) (F(i)-F(j))^2\\
%% \cvar_t(F) &= \sum_{i,j} \pi_t^{\bar\mu}(i) \pi_t^{\bar\mu}(j) (F(i)-F(j))^2
%\end{align*}
%is the same formula for both deterministic and random $F$.  
\end{remark}

\medskip

\begin{definition}\label{def:pathwisePI}
A Markov process is said to satisfy the {\em pathwise Poincar\'e
  inequality} with 
%respect to the filtration $\clZ$ and 
a constant $c$ if~\eqref{eq:cPI-pathwise}
holds and
\[
\beta_t \;\ge\; c\; >\; 0 \quad \sP^{\bar{\mu}}-\text{a.s.},\;\; t \ge 0
\]
\end{definition}

\medskip

The following proposition shows that conditional PI and pathwise PI
are equivalent.  Its proof appears in Appendix~\ref{apdx:pf-pathwise-PI-and-cPI}.  

\medskip

\begin{proposition}\label{prop:pathwise-PI-and-cPI}
The conditional PI holds with a constant $c$ if and only if the
pathwise PI holds with a constant $c$.  
\end{proposition}

\medskip

\begin{remark}\label{rem:reviewer_corrected}
For the conditional PI to hold, $\pi_t^{\bar{\mu}}$ must therefore satisfy the
standard form of the PI with a uniform constant $c$: 
\[
\sum_x  \pi_t^{\bar{\mu}}(x) \Gamma (f) (x) \; \geq \; c \;\sum_{x} \pi_t^{\bar{\mu}}(x)  |f(x) - \pi_t^{\bar{\mu}}(f)|^2 \quad \sP^{\bar{\mu}}-\text{a.s.}
\]
for all $f \in L^2(\bS)$ and for all $t\geq 0$.   This 
is a very stringent requirement that greatly limits the application of this
paper: It is easy to come up with examples where the standard form of
the PI holds for $\bar{\mu}$ but not for all $\pi_t^{\bar{\mu}}$.
The weaker form presented in~\Prop{prop:beta_ind_ineq},
specifically~\eqref{eq:PI_ind_ineq_filter_stronger}, is more general.
However, it is not clear how to obtain useful asymptotic bounds for $\frac{1}{T}\int_0^T
\beta_t \ud t$, beyond the result in Example~\ref{ex:min_row}.     
% In general, it is difficult to compute or bound the stochastic process
% defined by the righthand-side~\cite{pathiraja2020mckean}.  From
% \Prop{prop:pathwise-PI-and-cPI}, 
% conditional PI essentially requires that PI holds pathwise for all
% $\pi_t^{\bar{\mu}}$.  In general, this is a very stringent
% requirement.   
\end{remark}

%
%\begin{remark} For the conditional PI to hold, it is sufficient but
%  not necessary for the ratio defined by the righthand-side
%  of~\eqref{eq:cond_Rayleigh} to be bounded away from zero,
%  i.e., $\beta_t \ge c > 0$. This was the case for
%  Examples~\ref{ex:2-states},~\ref{ex:exdoeblin} and~\ref{ex:minsqrt},
%  but may be conservative.  The conditional PI only
%  requires that the ratio of averages 
%\[
%\frac{\den_t(F)}{\dvar_t(F)} = \frac{ \E^{\bar\mu} (\cen_t(F))}{ \E^{\bar\mu} (\cvar_t(F))} 
%\geq c > 0\quad \forall\; t\geq 0
%\]
%Of course the bound suggested
%in~\eqref{eq:PI_ind_ineq_filter_stronger} is even weaker.  It merely requires that $\int_0^\infty \beta_t \ud t =\infty$
%$\sP^{\bar{\mu}}$-a.s. for $\dvar_0({Y}_0)\to 0$.
%\end{remark}

\subsection{$Y_T$ as likelihood ratio}
In the remainder of this section $\pi^\mu$ and $\pi^{\bar\mu}$ are the
solutions of the Wonham filter~\eqref{eq:Kol}, with initialization
$\pi_0=\mu$ and $\pi_0=\bar{\mu}$, respectively.  
The likelihood ratio is the random function
\[
\gamma_T(x):=\frac{\pi_T^\mu(x)}{\pi_T^{\bar{\mu}}(x)}\quad \text{for} \;\; x\in\bS
\]  
% In contrast to \Sec{sec:stoch_stab}, the ratio now is a random
% function. 
The function is well-defined because $\pi_T^{\bar{\mu}}(x)>0$ for all
$x\in\bS$~\cite[Remark 3]{baxendale2004asymptotic}.  
 % {\color{red} Jin
 %  Kim.  We need to discuss why this is in the appropriate function
 %  space.  This is related to the mean and the variance being
 %  finite.}.
Its conditional mean and variance are
as follows:
\begin{align*}
& \text{mean}:\quad\quad \;\;\; \pi_T^{\bar\mu}(\gamma_T)   = \sum_x \pi_T^{\bar\mu}(x) \gamma_T(x)   =
  1\\
& \text{variance}:\quad\dvar_T({\gamma}_T)    = \E^{\bar{\mu}}(
                                         |\gamma_T(X_T)-1|^2) = \E^{\bar{\mu}}( \pi_T^\mu (\gamma_T) - 1)
\end{align*}
To derive the formula for the variance note
\begin{align*}
\dvar_T(\gamma_T) + 1 =
  \E^{\bar\mu}(|\gamma_T(X_T)|^2) &=
  \E^{\bar\mu}(\E^{\bar\mu}(|\gamma_T(X_T)|^2|\clZ_T)) \\ & =
\E^{\bar\mu}( \pi_T^{\bar\mu} (\gamma_T^2)) = \E^{\bar{\mu}}( \pi_T^\mu (\gamma_T))
\end{align*}
%Note that $\cvar_T(\gamma_T)$ is the Pearson's $\chi^2$-divergence.

% In order to use the inequality~\eqref{eq:PI_ind_ineq_filter}, it is
% natural to consider the dual optimal control problem with
% $Y_T=\gamma_T$.  Implications of this choice are described in the

Useful formulae for $\pi_T^\mu (\gamma_T)$ are obtained by 
setting $Y_T=\gamma_T$ in the dual optimal control problem. The results are
presented in the
following proposition whose
proof appears in the Appendix~\ref{appdx:filter_stab_prop}.  To state
the proposition, we need the exponential $\sP^{\bar\mu}$-martingale
\[
A_t := \exp\Big(\int_0^t D_s (
\ud Z_s - \pi_s^{\bar{\mu}}(h)\ud s) - \half \int_0^t D_s^2 \ud s
\Big),\quad t\geq 0
\]
where the difference
$D_t:=\pi_t^\mu(h)-\pi_t^{\bar{\mu}}(h)$. ($\{A_t:t\geq 0\}$ is
  in fact the change of measures between
%R-N derivative $\frac{\ud \sP^\mu_{\clZ_t}}{\ud     \sP^{\bar\mu}_{\clZ_t}}$ for
the conditional laws~\cite[Theorem 3.1]{clark1999relative}, but we do
  not make use of this interpretation here.)

\medskip

\begin{proposition} \label{prop:filter_stab_dual}
Consider the dual optimal control problem with
  $Y_T:=\gamma_T$.  Then
\begin{enumerate}
\item The optimal control ${U}^{\text{opt}}=0$ a.s., and the optimal
  trajectory $({Y},{V})$ is the 
  solution of the BSDE
\begin{align}\
- \ud {Y}_t(x)    & = \left((A {Y}_t)(x) + h (x) 
          {V}_t(x)\right)\ud t - {V}_t(x)\ud Z_t, 
%\; {Y}_T = \gamma_T
\nonumber \\
Y_T (x) &= \gamma_T(x) \;\; \forall \; x \in \bS \label{eq:opt_controlled_bsde_Uzero}
\end{align}
\item The process $\pi_t^{\bar\mu}({Y}_t) \equiv 1$ for
  all $0\leq t\leq T$.  Therefore
\[
\pi_T^{\bar\mu}(\gamma_T) = \bar\mu({Y}_0) = 1
\]
% and 
% \[
% \text{var}_T^{\bar{\mu}}(\bar{Y}_T) = \E^{\bar\mu}(|\gamma_T(X_T)-1|^2) 
% \]

\item The stochastic process $\{\pi_t^{\mu}({Y}_t): 0\leq t \leq
  T\}$ is a ${\sf P}^\mu$-mg. 
% \footnote{The proof would be easy if
    % this was ${\sf P}^{\bar\mu}$-mg.  See the next item in the list.}. 
  Consequently
\[
{\sf E}^{\mu} (\pi_T^{\mu}(\gamma_T)) % = {\sf E}^{\mu}
% (\pi_T^{\mu}(\gamma_T)) 
= \mu({Y}_0)
\] 

\item The stochastic process $\{A_t \pi_t^{\mu}({Y}_t): 0\leq t \leq
  T\}$ is a ${\sf P}^{\bar\mu}$-mg.  Consequently
\[
{\sf E}^{\bar\mu} (A_T \pi_T^{\mu}(\gamma_T)) % = {\sf E}^{\mu}
% (\pi_T^{\mu}(\gamma_T)) 
= \mu({Y}_0)
\] 
\end{enumerate}
\end{proposition}

\medskip

Now, using the inequality~\eqref{eq:PI_ind_ineq_filter} together with
the result from either part (3) or part (4),  it is a straightforward
calculation to derive the following forward inequality for the variance:
\begin{equation}\label{eq:RT_ineq}
R_T \; \dvar_T(\gamma_T) \leq e^{-cT}\dvar_0(\gamma_0)
\end{equation}
where $R_T$ is a factor obtained from either part~(3) or part~(4) of
the Proposition. The formulae for the two cases are:
\begin{itemize}
\item $R_T$ using part~(3):
\[
R_T = \left( \frac{\E^{{\mu}}( \pi_T^\mu ({\gamma}_T) -
    1)}{\E^{{\bar\mu}}( \pi_T^\mu ({\gamma}_T) - 1)} \right)^2
\]
\item $R_T$ using part~(4):
\[
R_T = \left( \frac{\E^{{\bar\mu}}( A_{T} (\pi_T^\mu ({\gamma}_T) -
    1))}{\E^{{\bar\mu}}( \pi_T^\mu ({\gamma}_T) - 1)} \right)^2
\]
\end{itemize}
% \[
% R_T = 
% \begin{cases}
%                     \left( \frac{\E^{{\mu}}( \pi_T^\mu ({\gamma}_T) -
%     1)}{\E^{{\bar\mu}}( \pi_T^\mu ({\gamma}_T) - 1)} \right)^2  & \text{using part (3)}  \\[10pt]
%                      \ \left( \frac{\E^{{\bar\mu}}( A_{T} (\pi_T^\mu ({\gamma}_T) -
%     1))}{\E^{{\bar\mu}}( \pi_T^\mu ({\gamma}_T) - 1)} \right)^2 & \text{using part (4)} 
%                  \end{cases} 
% \]
% where $\text{var}^{\bar{\mu}}(\gamma_0):=\bar{\mu}
%   ((\frac{\mu}{\bar\mu}-1)^2) = \mu(\gamma_0)-1$.  
The calculation for~\eqref{eq:RT_ineq} appears in the Appendix~\ref{appdx:filter_stab_prop} after the proof
  of the \Prop{prop:filter_stab_dual}.

\subsection{Filter stability}

We are interested in the difference
\begin{align*}
\pi_T^\mu(f) -  \pi_T^{\bar\mu}(f)  &= \pi_T^{\bar\mu}( (\gamma_T-1)
(f-\bar\mu(f))) \\
&= \E^{\bar\mu} ((\gamma_T(X_T)-1) (f(X_T)-\bar\mu(f))|\clZ_T)
\end{align*}
Taking expectations of the absolute value of both sides and using Cauchy-Schwarz
\begin{align*}
(\E^{\bar\mu}(|\pi_T^\mu(f) -  \pi_T^{\bar\mu}(f)|))^2 
\leq  \dvar_T(\gamma_T) \dvar_0(f)
% & =
%                                                      \E^{\bar\mu} \left(|\E^{\bar\mu}
%                                                      ((\gamma_T(X_T)-1)
%                                                      f(X_T)|\clZ_T)|\right)\\
% & \leq \E^{\bar\mu}\left( \E^{\bar\mu}
%                                                      (|(\gamma_T(X_T)-1)
%                                                      f(X_T)||\clZ_T)|\right)\\
% & = \E^{\bar\mu}(|(\gamma_T(X_T)-1)
%                                                      f(X_T)|)\\
% & \leq (\E^{\bar\mu}(|(\gamma_T(X_T)-1)|^2))^\half
% (\E^{\bar\mu}(|f(X_T)|^2))^\half
\end{align*}

Combined with the forward inequality~\eqref{eq:RT_ineq} for the
variance, we have proved the following main result for filter
stability.

\medskip

\begin{theorem}
Suppose the filter satisfies the conditional PI
with constant $c$.  Then   
\[
R_T (\E^{\bar\mu}(|\pi_T^\mu(f) -  \pi_T^{\bar\mu}(f)|))^2 
\leq e^{-cT}  \dvar_0(\gamma_0) \dvar_0(f)
\]
Consequently if $R_T \geq a^2 > 0$ then 
\[
(\E^{\bar\mu}(|\pi_T^\mu(f) -  \pi_T^{\bar\mu}(f)|))^2 
\leq \frac{1}{a^2} e^{-cT}  \dvar_0(\gamma_0) \dvar_0(f)
\]
\end{theorem}

\medskip

For equivalent measures, a conservative lower bound for $R_T$ is given
in the following proposition, whose proof appears in the
Appendix~\ref{appdx:prop:elem}.

\medskip

\begin{proposition}\label{prop:elem}
Suppose 
% the probability measures $\mu,\bar\mu\in\clP(\bS)$ are
% equivalent and let
$
a= \min_x \frac{\mu(x)}{\bar\mu(x)}>0
$.  
Then $R_T \geq a^2$. 

% Let $a^{-1} := \max_{i} \left\{ \left(\frac{\bar\mu(i)}{\mu(i)}\right)^2
% \right\}$. Then $R_T\geq a$.  Consequently, if the state space is
% finite and the measures $\mu, \bar\mu$ are equivalent then
% \[
% \text{var}_T^{\bar{\mu}}(\gamma_T) \to 0 \quad \text{as} \;\; T\to \infty
% \]
\end{proposition}

\medskip

\subsection{Forward variance inequality and filter stability}

Equation~\eqref{eq:RT_ineq} is the key variance inequality to obtain
the filter stability result.  Its counterpart for stochastic stability
is~\eqref{eq:RT_ineq_markov}.  For the ease of reader, these are
stated again:
\begin{align*}
\text{Markov process:}&\qquad \dvar_0(\gamma_T)  \leq e^{-cT} \dvar_0(\gamma_0)\\
\text{Wonham filter:}&\qquad R_T \; \dvar_T(\gamma_T) \leq e^{-cT}\dvar_0(\gamma_0)
\end{align*}  
The first of these inequalities is readily verified from a direct
calculation, which in fact is a standard proof to prove
stochastic stability using the PI (see~\cite[Theorem 4.2.5]{bakry2013analysis}):
\[
\frac{\ud }{\ud t}\dvar_0(\gamma_t) = - \den_0(\gamma_t) 
\]
and using PI, the variance inequality for the Markov process follows.  In contrast, a
duality based proof of the variance inequality, described in
Appendix~\ref{sec:stoch_stab}, is more involved.  

For the filter, the variance inequality~\eqref{eq:RT_ineq} is new.
One may ask whether~\eqref{eq:RT_ineq} can also be derived
more directly (without the use of duality)?  A direct calculation
shows that
\[
\ud \cvar_t(\gamma_t) = - \cen_t(\gamma_t)\ud t + C_t \ud t +
(\sP^\mu\text{-mg. increment})
\]
where the coefficient 
\[
C_t:=\pi_t^{\bar{\mu}}\big(\gamma_t^2(h-\pi_t^\mu(h))\big)(\pi_t^{\bar{\mu}}
  (h)-\pi_t^\mu(h))
\]
It is not clear how the equation can be simplified to
obtain~\eqref{eq:RT_ineq}.  It may be possible to derive~\eqref{eq:RT_ineq} through a clever choice of an integrating factor.  However,
we have not yet been successful in this endeavor.  
Because the BSDE~\eqref{eq:opt-cont-euclide-b} and the Zakai
equation are dual (see~\cite[Prop. 1]{kim2019observability}), it
will be helpful to relate this work to the study of the Lyapunov
exponents of the Zakai equation~\cite{atar1997lyapunov,atar1997exponential,Atar1999}.

\section{Acknowledgement}

Comments from an anonymous reviewer helped greatly improve this paper.
In particular, the reviewer corrected an error in
a prior version by pointing that the conditional PI and pathwise
PI are equivalent (\Prop{prop:pathwise-PI-and-cPI}).  
The implications of this equivalence are discussed in Remark~\ref{rem:reviewer_corrected}.    

\bibliographystyle{IEEEtran}
\bibliography{duality,backward_sde,filter-stability-observability,jin_papers}

\appendix

\section{Appendix}

\subsection{Proof of \Prop{prop:CPI-examples}}
\label{appdx:CPI-examples-pf}

Because the formulae are of independent interest, their proofs are
included as part of two examples which are described below.  Together
with Examples~\ref{ex:2-states} and~\ref{ex:min_row}, they comprise
the four examples in the paper where the conditional PI holds. 
\begin{example}\label{ex:exdoeblin}
A Markov process is {\em Doeblin} if there exist a state $j^*\in \bS$
such that $A(x,j^*) > c$ for all $x\in \bS \setminus \{j^*\}$.
Then 
\begin{align*}
\den_T(F) & = 
  \E^{\bar\mu} \Big( \sum_{j \in \mathbb{S}} A(X_T,j) (F(X_T) - F(j))^2 \Big)\\
  &\geq c\,
  \E^{\bar\mu} ( (F(X_T) - F(j^*))^2 ) 
%\\ 
%& \geq \epsilon \; \E^{\bar\mu} ( (F(X_T) - \pi_T^{\bar\mu}(F) )^2 )  
\geq c \,  \dvar_T(F)
\end{align*}
Therefore, the conditional PI holds for Doeblin chains.  This admits a
straightforward generalization that gives the first of the two
formulae in \Prop{prop:CPI-examples}:  
\begin{align*}
\den_T(F) & \geq  
  \E^{\bar\mu} \Big( \sum_{j \in \mathbb{S}} \min_{i:\;i\neq j} A(i,j)
            \; (F(X_T) - F(j))^2 \Big)\\
  &\geq \Big( \sum_{j \in \mathbb{S}} \min_{i:\;i\neq j} A(i,j) \Big) \dvar_T(F)
%\\ 
%& \geq \epsilon \; \E^{\bar\mu} ( (F(X_T) - \pi_T^{\bar\mu}(F) )^2 )  
%\geq c \,  \dvar_T(F)
\end{align*}
\end{example}

\medskip

\begin{example}\label{ex:minsqrt}
% Because  algebraic mean dominates the geometric mean
% \begin{align*}
% \den_T(F) & = 
%   \E^{\bar\mu} \Big( \sum_{j} \frac{1}{2}
%             (A(X_T,j)+A(j,X_T)) (F(X_T) - F(j))^2 \Big)\\
% &\geq \E^{\bar\mu} (\sum_{j} \sqrt{A(X_T,j) A(j,X_T)} (F(X_T) - F(j))^2
% \end{align*}
Using the tower property of conditional expectation
\begin{align*}
\den_T(F) &= \E^{\bar{\mu}} \Big(\sum_{i,j\in \bS} \pi_T^{\bar{\mu}}(i)A(i,j)(F(i)-F(j))^2\Big)\\
\dvar_T(F) &= \E^{\bar{\mu}} \Big(\sum_{i,j\in\bS} \pi_T^{\bar{\mu}}(i)\pi_T^{\bar{\mu}}(j) (F(i)-F(j))^2\Big)
\end{align*}
Because algebraic mean dominates the geometric mean
\begin{align*}
\sum_{i,j\in \bS} &\pi_T^{\bar{\mu}}(i)A(i,j)(F(i)-F(j))^2\\
&=\sum_{i,j\in \bS} \frac{1}{2}\big(\pi_T^{\bar{\mu}}(i)A(i,j)+\pi_T^{\bar{\mu}}(j)A(j,i)\big)(F(i)-F(j))^2\\
&\ge\sum_{i,j\in \bS} \sqrt{\pi_T^{\bar{\mu}}(i)\pi_T^{\bar{\mu}}(j)}\sqrt{A(i,j)A(j,i)}(F(i)-F(j))^2\\
&\ge c\sum_{i,j\in\bS} \pi_T^{\bar{\mu}}(i)\pi_T^{\bar{\mu}}(j) (F(i)-F(j))^2
\end{align*}
where we used the fact that $\sqrt{x} \geq x$ for $0\leq x\leq 1$.  Taking
expectations of both sides gives the second of the two fromulae in \Prop{prop:CPI-examples}. 
\end{example}

\subsection{A counterexample for conditional PI}
\label{apdx:counterexample}
\begin{example}
This is the famous counterexample of the filtering
theory~\cite[pp. 9-10]{delyon1988lyapunov}.  The state-space
$\bS=\{1,2,3,4\}$ and the rate matrix 
\[
A = \begin{pmatrix}
-1 & 1 & 0 & 0\\
0 & -1 & 1 & 0\\
0 & 0 & -1 & 1\\
1 & 0 & 0 & -1
\end{pmatrix}
\]
whose unique invariant measure $\bar{\mu} =
[\frac{1}{4},\frac{1}{4},\frac{1}{4},\frac{1}{4}]$.  The standard PI
holds with a constant $c= 2$.

Consider a sigma-algebra $\clG = \sigma([X_T \in\{1,3\}])$ along with
a $\clG$-measurable function:
\[
F(\cdot) = \begin{cases}
                    \begin{pmatrix} 1 & 1 & -1 & -1\end{pmatrix} & \text{if}  \; \; X_T \in\{1,3\}  \\[4pt]
                     \begin{pmatrix} -1 & 1 & 1 & -1 \end{pmatrix} & \text{if}  \; \;  X_T
                     \in\{2,4\} 
                 \end{cases} 
\]
Then the conditional distribution
\[
\pi_T^{\bar\mu}(\cdot) = \E^{\bar{\mu}} (\ones_{[X_T=\cdot]}\mid \clG)
= \begin{cases}
                    \begin{pmatrix} \frac{1}{2} & 0 & \frac{1}{2} & 0\end{pmatrix} & \text{if}  \; \; X_T \in\{1,3\}  \\[4pt]
                     \begin{pmatrix} 0 & \frac{1}{2} & 0 &
                       \frac{1}{2} \end{pmatrix} & \text{if}  \; \;  X_T
                     \in\{2,4\} 
                 \end{cases} 
\]
The conditional mean $\pi_T^{\bar\mu}(F)=0$, the conditional
variance $\pi_T^{\bar\mu}(F^2)=1$, and therefore the variance
$\dvar_T(F)=1$.  On the other hand, the energy
$
\den_T(F) = 0
$.  
  Therefore,
the conditional PI does not hold for this example.  
\end{example}

% \begin{remark}
% This is the famous counterexample on stochastic stability, first noted by~\cite{delyon1988lyapunov}. In section 3 of \cite{baxendale2004asymptotic}, it is pointed out that the main conjecture in~\cite{kunita1971asymptotic} fails for this example with noiseless observation, and the filter is not stable even though the state process is ergodic.
% Additional assumption on the observation is provided in~\cite[Assumption III.2]{van2010nonlinear}.
% \end{remark}

\subsection{Stochastic stability}\label{sec:stoch_stab}
 
% In this section, a proof for the stochastic stability of Markov
% processes is described based upon assuming the standard PI and dual
% formulation. This lies parallel with the discussion on
% Sec.~\ref{sec:filt_stab}.

In this section, we specialize the results of Sec.~\ref{sec:filt_stab}
to the problem of stochastic stability of Markov processes. 

\newP{Dual process}
Introduce a deterministic dual process $y=\{y_t\in L^2(\bS):0 \leq t \leq T\}$
as a solution of the backward ode
\[
- \frac{\ud }{\ud t} y_t = A y_t, \quad y_T  \; \text{fixed (given)}
\]
Then a standard application of It\^o product formula gives
\begin{equation*}\label{eq:eq_f_st}
y_T(X_T) = y_0(X_0) + \int_0^T \sum_{x\in\bS} y_t(x)\ud N_t(x)
%\langle y_t,\ud N_t\rangle
\end{equation*}
where 
%$\langle y_t,\ud N_t\rangle=\sum_{x\in\bS} y_t(x)\ud N_t(x)$ and  
$\{N_t:t\geq 0\}$ is a martingale.  
% %$N_t := X_t - \int_0^t A^\top X_s \ud s$
% $\{N_t:t\geq 0\}$ is a martingale (mg) and $\langle y_t,\ud
% N_t\rangle=BLAH$. 
% \marginpar{Fix BLAH}
% $\{N_t:t\geq 0\}$ is a martingale characterization associated with $X$ defined for each values of $x\in\bS$ by
% \[
% N_t (x) := \ones_{X_t}(x) -\ones_{X_0}(x) - \int_0^t A(X_s,x) \ud s
% \]
Taking an expectation
\begin{equation*}\label{eq:eq_f_dt}
\E^{\bar\mu} (y_T (X_T)) = \bar\mu(y_0)
\end{equation*}
and the equation for the variance is
%Upon subtracting~\eqref{eq:eq_f_dt} from~\eqref{eq:eq_f_st}
% \[
% y_0(X_0) - \bar\mu(y_0) + \int_0^T \langle y_t,\ud N_t\rangle = y_T (X_T) - \E^{\bar\mu} (y_T (X_T))
% \]
% and upon squaring and taking expectations
\[
\dvar_0(y_0) + \int_0^T \den_0(y_t) \ud t  = \dvar_0(y_T)
\]

% \begin{equation*}\label{eq:eq_f_var}
% \E^\mu|y_0(X_0) - \mu(y_0)|^2  +
% \int_0^T \E^\mu (\Gamma (y_t) (X_t))  \ud t = \E^\mu |y_T (X_T) - \E^\mu (y_T (X_T))|^2 
% \end{equation*}
% where \begin{equation}\label{eq:Q-explicit}
% Q(e_i) := \sum_{j=1}^dA_{ij}(e_j-e_i)(e_j-e_i)^\top
% \end{equation}
% is the carre du champs operator.  
%\end{enumerate}

% Taking $\mu=\bar{\mu}$, the equation for the variance reads
% \[
% \text{var}^{\bar{\mu}}(y_0) + \int_0^T {\cal E}^{\bar{\mu}}(y_t) \ud t  = \text{var}^{\bar{\mu}}(y_T)
% \]
Using the PI it follows
\begin{equation}\label{eq:PI_ind_ineq_markov}
\dvar_0(y_0) \leq e^{-cT} \dvar_0(y_T)
\end{equation}
This is counterpart of~\eqref{eq:PI_ind_ineq_filter}.  
This inequality is the only place in the stability proof where the PI
is used.  % For this reason, Eq.~\eqref{eq:PI_ind_ineq_markov} is
% referred to as the PI-induced inequality. 

% This inequality is the counterpart
% of~\eqref{eq:PI_ind_ineq_filter}. 

\newP{$\mathbf{y_T}$ as likelihood ratio}
In the remainder of this section $\pi^\mu$ and $\pi^{\bar\mu}$ are
the solutions of the Kolmogorov's equation~\eqref{eq:Kol}, with initialization
$\pi_0=\mu$ and $\pi_0=\bar{\mu}$, respectively.  

Consider the likelihood ratio
$
\gamma_T(x):=\frac{\pi_T^\mu(x)}{\bar{\mu}(x)}\; \text{for} \;\; x\in\bS
$.  
% (this means
% $\gamma_T(x)=\frac{\pi_T^\mu(x)}{\bar{\mu(x)}}$ for all $x\in\bS$).
The ratio is well-defined because $\bar{\mu}(x)>0$ for all
$x\in\bS$.  The mean and variance of $\gamma_T$ are as follows:
\begin{align*}  
& \text{mean}:\quad\quad \;\;\; \bar\mu(\gamma_T)   = \sum_x
  \bar\mu(x) \gamma_T(x) = 1 \\
& \text{variance}:\quad
\dvar_0(\gamma_T)  = {\bar \mu} (\gamma_T^2) - 1 =
                                    \pi_T^\mu (\gamma_T) - 1
\end{align*}

%The counterpart of \Prop{prop:prop1} is as follows:
The dual process is used to express $\pi_T^\mu (\gamma_T)$ in terms of
the initial measure $\mu$.  The counterpart
of~\Prop{prop:filter_stab_dual} is as follows:

% The proof of the following proposition appears in the
% Appendix~\ref{appdx:stoch_stab}. 
\medskip

\begin{proposition}\label{prop:stoch_stab_dual}
Consider the dual process with terminal condition $y_T:=\gamma_T$.  Then
\begin{align*}
\bar\mu(\gamma_T) &= \bar\mu(y_0) ,\quad 
\pi_T^\mu (\gamma_T)  = \mu(y_0)
\end{align*}
\end{proposition}

\medskip

\begin{proof}
Let $\{\pi_t:0\leq t \leq T\}$ be the solution of the Kolmogorov's forward equation then 
\[
\pi_T(y_T) = \pi_0(y_0)
\]
The two formulae 
follow by using $\pi_0=\bar{\mu}$ and $\pi_0=\mu$.  % This
% concludes the proof of~\Prop{prop:stoch_stab_dual}. 
\end{proof}
\medskip

Using the inequality~\eqref{eq:PI_ind_ineq_markov} together with
the result of the proposition, the following counterpart
of~\eqref{eq:RT_ineq} is obtained:
\begin{equation}\label{eq:RT_ineq_markov}
\dvar_0(\gamma_T)  \leq e^{-cT} \dvar_0(\gamma_0)
\end{equation}
The calculation is included at the end of this section.  % This inequality is the

\newP{Stochastic stability}
%The stability result now follows readily.  
Using Cauchy-Schwarz, the difference squared 
\[
|\pi_T^\mu(f) -  \bar{\mu}(f)|^2 = |\bar{\mu}((\gamma_T-1)(f-\bar\mu(f)))|^2
\leq \dvar_0(\gamma_T) \dvar_0(f)
\]  
% Using Cauchy-Schwarz,
% \[
% |\bar{\mu}((\gamma_T-1)f)|^2 \leq \text{var}^{\bar{\mu}}(\gamma_T)\,
% \bar{\mu} (f^2) 
% \]
Combining this with~\eqref{eq:RT_ineq_markov}, we have the following result on stochastic stability:

\medskip

\begin{theorem}
Suppose $X$ is a Markov process with an everywhere positive invariant
measure $\bar\mu$ whose generator satisfies the PI with constant $c$. Then
%  Let $\mu\ll\bar\mu$ with $\bar{\mu}
% (\frac{\mu^2}{\bar\mu^2})$ finite and suppose $f$ is a given function with
% $\bar\mu(f^2)<\infty$ then
\[
|\pi_T^\mu(f)  -  \bar{\mu}(f)|^2 \leq  e^{-cT}\dvar_0(\gamma_0)
\dvar_0(f)  
\]
\end{theorem}

% \subsubsection{Proof of \Prop{prop:stoch_stab_dual} and Eq.~\eqref{eq:RT_ineq_markov}}
% \label{appdx:stoch_stab}

% We next prove the
% inequality in Eq.~\eqref{eq:RT_ineq_markov}.  

\medskip

\newP{Calculation for~\eqref{eq:RT_ineq_markov}}
Using the second formula from ~\Prop{prop:stoch_stab_dual}
%the variance of the likelihood ratio
\[
\dvar_0(\gamma_T) = \pi_T^\mu (\gamma_T) - 1 =  \mu(y_0)- 1 = \bar{\mu}\Big(\big(\frac{\mu}{\bar\mu}-1\big)(y_0-1)\Big) 
\]
and using Cauchy-Schwarz
\[
(\dvar_0(\gamma_T))^2 \leq \bar{\mu}
  \Big(\big(\frac{\mu}{\bar\mu}-1\big)^2\Big) \; \bar{\mu}
                               \big((y_0 - 1)^2\big) = 
  \dvar_0(\gamma_0) \; \dvar_0(y_0)  
\]

Finally, use the inequality~\eqref{eq:PI_ind_ineq_markov} for the
dual process:
\begin{align*}
\dvar_0(y_0) & \leq e^{-cT} \dvar_0(\gamma_T) \\
\therefore \;\; 
(\dvar_0(\gamma_T))^2 & \leq
\dvar_0(\gamma_0)
 \; \dvar_0(y_0)\leq  e^{-cT} \dvar_0(\gamma_0) \dvar_0(\gamma_T)
 \\
\therefore \quad\;\;\; \dvar_0(\gamma_T) & \leq e^{-cT} \dvar_0(\gamma_0)
\end{align*}
The proves the inequality in~\eqref{eq:RT_ineq_markov}.

\subsection{Proof of \Prop{prop:beta_ind_ineq}}\label{appdx:beta_ind_ineq_pf}

Using the definitions for $M$ and $\beta$
\begin{align*}
M_t &  = \cvar_t(Y_t) - \int_0^t \big( \cen_s(Y_s) +
\sum_x \pi_s^{\bar\mu} (x) |U_s^{\text{opt}}+V_s(x)|^2_R \big) \ud s\\
& \leq \cvar_t(Y_t) - \int_0^t \cen_s(Y_s) \ud s\;\; \leq \;\;
  \cvar_t(Y_t) - \int_0^t \beta_s \cvar_s(Y_s) \ud s 
\end{align*}
%and $M_0 = \cvar_0(Y_0) = \dvar_0(Y_0)$.  
That is
\begin{equation}\label{eq:Mtineq}
M_t \leq \cvar_t(Y_t) - \int_0^t \beta_s \cvar_s(Y_s) \ud s \quad
\sP^{\bar\mu}-\text{a.s.}\quad \text{for}\;\;0\leq t\leq T
\end{equation}
% The inequalities hold
% $\sP^{\bar\mu}$-a.s for $0\leq t\leq T$ above and also everywhere in this section.  
Set $\Phi_t = e^{-\int_0^t \beta_s \ud s}$ and multiply both sides by
the integrating factor
$\beta_t\Phi_t$ to obtain
\[
\beta_t \Phi_t M_t  \leq \beta_t \Phi_t\Big(\cvar_t(Y_t) - \int_0^t \beta_s \cvar_s(Y_s)\ud s\Big) 
\]
and upon integrating from $0$ to $T$
\begin{align*}
\int_0^T \beta_t \Phi_t M_t \ud t \leq \Phi_T\int_0^T \beta_t \cvar_t(Y_t)\ud t 
\end{align*}
% This is because
% \[
% \beta_t \Phi_t\Big(\cvar_t(Y_t) - \int_0^t \beta_s \cvar_s(Y_s)\ud s\Big) = \frac{\ud}{\ud t} \Big(\Phi_t \int_0^t \beta_s \cvar_s(Y_s)\ud s\Big)
% \]
Because $\ud \Phi_t = - \beta_t \Phi_t \ud t$ 
% \begin{align*}
% \int_0^T \beta_t \Phi_t (\cvar_0(Y_0)+M_t)\ud t
%  &= -\Phi_T(\cvar_0(Y_0)+M_T) \\
%  &\quad +\cvar_0(Y_0) +  \int_0^T \Phi_t \ud M_t
% \end{align*}
\begin{align*}
-\Phi_T M_T +  M_0 +  \int_0^T \Phi_t \ud M_t \leq \Phi_T\int_0^T \beta_t \cvar_t(Y_t)\ud t 
% \Phi_T\Big(\int_0^T \beta_s \cvar_s(Y_s)\ud s + \cvar_0&(Y_0)+M_T\Big)\\
% &\ge\;\; \cvar_0(Y_0) +  \int_0^T \Phi_t \ud M_t
\end{align*}
and because $M_0 = \cvar_0(Y_0) = \dvar_0(Y_0)$
\begin{align*}
\dvar_0(Y_0) +  \int_0^T \Phi_t \ud M_t  & \leq \Phi_T\Big(\int_0^T
\beta_s \cvar_s(Y_s)\ud s + M_T\Big) \\ & \leq \Phi_T \cvar_T(Y_T)
\end{align*}
where~\eqref{eq:Mtineq} is used (with $t=T$) to obtain the final
inequality.  Take expectations of both sides to
obtain~\eqref{eq:PI_ind_ineq_filter_stronger}.  
% The left-hand side is bounded due to~\eqref{eq:beta_induced_ineq} evaluated at $T$, so
% \[
% \Phi_T\cvar_T(Y_T) \ge \cvar_0(Y_0) + \int_0^T \Phi_t \ud M_t
% \]
% Upon taking expectation, the claim is proved.

\subsection{Proof of Prop.~\ref{prop:pathwise-PI-and-cPI}}\label{apdx:pf-pathwise-PI-and-cPI}

Taking an expectation of both sides of~\eqref{eq:cPI-pathwise}, using
$\beta_t>c$, one obtains the conditional PI.  Therefore, pathwise PI
with a constant $c$ implies conditional PI with the constant $c$.
Conversely, consider a random function $F=1_B f$ where the set $B\in\clZ_T$
and $f\in L^2(\bS)$.  Using the property of the conditional
expectation $\cen_T (1_B  f) = 1_B\cen_T (f)$ and $\cvar_T(1_B f)= 1_B\cvar_T(f)$, and therefore conditional PI implies
\[
\E^{\bar\mu}\big(1_B\cen_T(f)\big) \; \ge \; c\,\E^{\bar\mu}\big(1_B\cvar_T(f)\big)
\]       
Since $B$ is arbitrary, pathwise PI follows.  

% If part is simple, because $F(\omega) \in L^2$, so the inequality holds point-wise.

% For only if part, assume the conditional PI. If $A
% \in\clZ_T$, then $1_Af\in L^2_{\clZ_T}$ for any $f\in L^2$. Observe that $\cen_T (1_Af) = 1_A\cen_T (f)$ and $\cvar_T(1_Af)= 1_A\cvar_T(f)$, and therefore
% \[
% \E^{\bar\mu}\big(1_A\cen_T(f)\big) \ge c\E^{\bar\mu}\big(1_A\cvar_T(f)\big)
% \]
% Set $A = [{\cal E}_t^{\bar\mu} (f) < c {\cal V}_t^{\bar\mu}(f)]$, then $\sP(A) = 0$ or the assumption is violated. %\qed

\subsection{Proof of Prop.~\ref{prop:filter_stab_dual} and~\eqref{eq:RT_ineq}}
\label{appdx:filter_stab_prop}

%\begin{proof}{of Prop.~\ref{prop:filter_stab_dual}}
\begin{enumerate}
\item With $Y_T=\gamma_T$, the conditional expectation
\begin{align*}
\pi_T^{\bar\mu}(Y_T)
= \sum_{x\in\bS} \pi_T^{\bar{\mu}}(x) \gamma_T(x)
= \sum_{x\in\bS}\pi_T^{\mu}(x) = 1
\end{align*}
Therefore $U_t^{\text{opt}}$ satisfies 
\[
1 = \bar{\mu}({Y}_0) - \int_0^T U_t^{\text{opt}} \ud Z_t \quad \sP^{\bar\mu}-\text{a.s.}
\]
If $Z$ was a w.p. then it follows from the representation
theorem~\cite[Theorem 5.18]{le2016brownian} that $U_t^{\text{opt}}
\equiv 0$ a.s. for $0\leq t\leq T$ and $\bar{\mu}({Y}_0) =1$. 
Moreover, the representation is unique in $L^2_\clZ([0,T])$.  We can
not apply the representation theorem directly because $Z$ is not a
$\sP^{\bar\mu}$-w.p.  However, it is w.p. under the Girsanov change of
measure~\cite[p. 85]{xiong2008introduction}.  Since the two measures
are equivalent, the representation also holds for $\sP^{\bar\mu}$.

\item 
In~\cite[Theorem~2]{kim2019duality}, it is proved that along the
optimal trajectory
\[
\pi_t^{\bar\mu}(Y_t)  = \bar{\mu}({Y}_0) - \int_0^t U_s^{\text{opt}}
\ud Z_s,\quad \text{for}\; \; 0\leq t \leq T
\]
Therefore, $\pi_t^{\bar\mu}(Y_t)  = \bar{\mu}({Y}_0) = 1$.  

\item % In this part and the next $\pi^\mu=\{\pi^\mu_t:0\leq t\leq T\}$ is a solution of the Wonham
  % filter with initialization $\pi_0=\mu$.  
  Using the equation of the
  Wonham filter~\eqref{eq:Wonham} for $\pi^\mu$ and the 
  BSDE~\eqref{eq:opt_controlled_bsde_Uzero} for $Y$, a direct calculation
  (included in Appendix~\ref{appdx:ofappdx})) shows 
\begin{equation}\label{eq:pitmuyt}
\ud (\pi_t^{\mu}({Y}_t)) = C_t  \big(\ud Z_t - \pi_t^\mu(h)\ud
t\big),\quad 0\leq t\leq T
\end{equation}
where the coefficient
\[C_t=\pi_t^\mu(h {Y}_t)-\pi_t^\mu({Y}_t)\pi_t^\mu(h)
+ \pi_t^\mu ( V_t) ,\quad 0\leq t\leq T
\] 

\item This part also follows from a direct calculation (included in Appendix~\ref{appdx:ofappdx})) to show that
\begin{align}\label{eq:atpitmuyt}
\ud \big(A_t \pi_t^\mu({Y}_t)\big) =
A_t (D_t\pi_t^\mu(Y_t)+ C_t) (\ud Z_t -
                                         \pi_t^{\bar\mu}(h)\ud t) 
\end{align}
where $C_t$ is as defined in the proof of part 3.   This 
completes the proof of the four parts of the Proposition.  
\end{enumerate}
%\end{proof}

\medskip

We next derive the inequality in~\eqref{eq:RT_ineq}.  Two
derivations are provided starting from results in part 3 and part 4:

%\medskip

\newP{Derivation of~\eqref{eq:RT_ineq} from part 3 of Prop.~\ref{prop:filter_stab_dual}} Since
$\{\pi_t^{\mu}({Y}_t):0\leq t\leq T\}$ is $\sP^\mu$-mg
\[
\E^{{\mu}}( \pi_T^\mu ({\gamma}_T) - 1) =
\mu(({Y}_0-1)) = \bar{\mu} ((\frac{\mu}{\bar\mu}-1)({Y}_0-1)) 
\]
and therefore using Cauchy-Schwarz
\[
|\E^{{\mu}}( \pi_T^\mu ({\gamma}_T) - 1)|^2 \leq \dvar_0(\gamma_0) \; \dvar_0({Y}_0)  
\]
%where $\text{var}^{\bar{\mu}}(\gamma_0)=\bar{\mu}
%((\frac{\mu}{\bar\mu}-1)^2)$.  

Now we use the inequality~\eqref{eq:PI_ind_ineq_filter} for the dual process
\begin{align*}
\dvar_0({Y}_0) & \leq e^{-cT}
                                      \dvar_T(\gamma_T)\\
\therefore, \quad |\E^{{\mu}}( \pi_T^\mu ({\gamma}_T) - 1)|^2 &
 % \leq \text{var}_0^{\bar{\mu}}(\gamma_0) \;
% \text{var}_0^{\bar{\mu}}(\bar{Y}_0)
\leq e^{-cT}\dvar_0(\gamma_0)  \dvar_T(\gamma_T)
\end{align*}
which is expressed as
\[
R_T \; \dvar_T(\gamma_T) \leq e^{-cT}\dvar_0(\gamma_0)
\]
with
\[
R_T:= \left( \frac{\E^{{\mu}}( \pi_T^\mu ({\gamma}_T) -
    1)}{\E^{{\bar\mu}}( \pi_T^\mu ({\gamma}_T) - 1)} \right)^2
\]

\newP{Derivation of~\eqref{eq:RT_ineq} from part 4 of Prop.~\ref{prop:filter_stab_dual}} Since
$\{A_t\pi_t^{\mu}({Y}_t):0\leq t\leq T\}$ is $\sP^{\bar\mu}$-mg
\[
\E^{{\bar\mu}}( A_T(\pi_T^\mu ({\gamma}_T) - 1)) =
\mu(({Y}_0-1)) = {\mu} ((\frac{\mu}{\bar\mu}-1)({Y}_0-1)) 
\]
and therefore using Cauchy-Schwarz
\[
|\E^{{\bar\mu}}( A_T  (\pi_T^\mu ({\gamma}_T) - 1))|^2 \leq \dvar_0(\gamma_0) \; \dvar_0({Y}_0)  
\]
Now we use the inequality~\eqref{eq:PI_ind_ineq_filter} for the dual process
\begin{align*}
\dvar_0({Y}_0) & \leq e^{-cT}
                                      \dvar_T(\gamma_T)\\
\therefore, \quad |\E^{{\bar\mu}}( A_T  (\pi_T^\mu ({\gamma}_T) - 1))|^2  & 
% \leq \text{var}_0^{\bar{\mu}}(\gamma_0) \;
% \text{var}_0^{\bar{\mu}}(\bar{Y}_0)
\leq e^{-cT}\dvar_0(\gamma_0) \dvar_T(\gamma_T)
\end{align*}
which is expressed as
\[
R_T \; \dvar_T(\gamma_T) \leq e^{-cT}\dvar_0(\gamma_0)
\]
with
\[
R_T:= \left( \frac{\E^{{\bar\mu}}( A_{T} (\pi_T^\mu ({\gamma}_T) -
    1))}{\E^{{\bar\mu}}( \pi_T^\mu ({\gamma}_T) - 1)} \right)^2
\]

\subsection{Proof of Prop.~\ref{prop:elem}}
\label{appdx:prop:elem}

Consider the random variable $S_T:=\pi_T^\mu(\gamma_T)-1$.  The ratio
$R_T = \left( \frac{\E^{{\mu}}(S_T)}{\E^{{\bar\mu}}( S_T)}
\right)^2$.  Using the law of total probability
\begin{align*}
\E^{{\mu}}( S_T )  = \sum_x \E^{{\delta_x}}( S_T ) \mu(x),\quad 
\E^{{\bar\mu}}( S_T ) = \sum_x \E^{{\delta_x}}( S_T ) \bar\mu(x)
\end{align*}
Therefore,
\[
\E^{{\mu}}( S_T ) = \sum_x \E^{{\delta_x}}( S_T )
\frac{\mu(x)}{\bar\mu(x)} \bar\mu(x) \geq a \E^{{\bar\mu}}( S_T ) 
\]
The bound follows.

\subsection{Additional calculations}
\label{appdx:ofappdx}

\newP{Calculation for~\eqref{eq:pitmuyt}} For this and the next
calculation, $\pi^\mu$ is a solution of the Wonham filter~\eqref{eq:Wonham} and $Y$ is a
solution of the BSDE~\eqref{eq:opt_controlled_bsde_Uzero}. Using the It\^{o}-Wentzell formula for measures~\cite[Theorem 1.1]{krylov2011ito}:
\begin{align*}
\ud \pi_t^\mu(Y_t) 
% &= \frac{\ud\sigma_t^\mu(Y_t)}{\sigma_t^\mu(1)} - \frac{\sigma_t^\mu(Y_t)}{(\sigma_t^\mu(1))^2} \ud \sigma_t^\mu(1) \\
% &\quad+ \frac{\sigma_t^\mu(Y_t)}{(\sigma_t^\mu(1))^3} \ud \sigma_t^\mu(1)\ud \sigma_t^\mu(1) - \frac{\ud \sigma_t^\mu(Y_t)\ud \sigma_t^\mu(1)}{(\sigma_t^\mu(1))^2}\\
&= \pi_t^\mu\big(h Y_t + V_t\big)\ud Z_t - \pi_t^\mu(Y_t)\pi_t^\mu(h) \ud Z_t \\
&\quad + \pi_t^\mu(Y_t)\pi_t^\mu(h) \pi_t^\mu(h)\ud t -
  \pi_t^\mu\big(h Y_t + V_t\big)\pi_t^\mu(h)\ud t\\
&= \big(\pi_t^\mu\!(h Y_t)-\pi_t^\mu\!(Y_t)\pi_t^\mu\!(h) + \pi_t^\mu\!(V_t)\big)(\ud Z_t - \pi_t^\mu\!(h)\ud t) \\
&= C_t \big(\ud Z_t - \pi_t^\mu(h)\ud t\big)
\end{align*}

\newP{Calculation for~\eqref{eq:atpitmuyt}}
The exponential martingale $A_t$ is the Dol\'{e}ans exponential of $\int_0^t D_s(\ud Z_s - \pi_s^{\bar{\mu}}(h)\ud s)$, and therefore its differential form is given by
\[
\ud A_t = A_tD_t(\ud Z_t - \pi_t^{\bar{\mu}}(h)\ud t)
\]
Applying It\^{o}'s lemma together with~\eqref{eq:pitmuyt} yields
\begin{align*}
\ud \big(A_t \pi_t^\mu(Y_t)\big) &= A_t D_t\pi_t^\mu(Y_t)(\ud Z_t - \pi_t^{\bar{\mu}}(h)\ud t) \\
&\quad + A_tC_t\big(\ud Z_t - \pi_t^\mu(h)\ud t\big) + A_tC_tD_t \ud t\\
&= A_t\big(D_t\pi_t^\mu(Y_t) + C_t\big)(\ud Z_t - \pi_t^{\bar{\mu}}(h)\ud t)
\end{align*}

\end{document}